\input amstex
\documentstyle{amsppt}
\document
\magnification=1200
\NoBlackBoxes
\nologo
\pageheight{18cm}

%\hfill{\it file Amywork/segre.tex, Nov. 30, 2004}

\bigskip

\centerline{\bf MANIFOLDS WITH MULTIPLICATION} 

\smallskip

\centerline{\bf ON THE TANGENT SHEAF\footnotemark1}
\footnotetext{Talk at the Conference dedicated
to the memory of B.~Segre, Inst. Mat. Guido Castelnuovo,
Rome, June 2004.}

\medskip

\centerline{\bf Yuri I. Manin}

\medskip

\centerline{\it Max--Planck--Institut f\"ur Mathematik, Bonn, Germany,}

\centerline{\it and Northwestern University, Evanston, USA}

\bigskip

{\bf Abstract.} This talk reviews the current state
of the theory of $F$--(super)manifolds $(M,\circ )$,
first defined in [HeMa] and further developed in [He],
[Ma2], [Me1]. Here $\circ$ is an $\Cal{O}_M$--bilinear multiplication
on the tangent sheaf $\Cal{T}_M$, satisfying an
integrability condition. $F$--manifolds and compatible
flat structures on them furnish a useful weakening
of Dubrovin's Frobenius structure which naturally arises in the 
quantum $K$--theory, theory of extended moduli spaces,
and unfolding spaces of singularities.

\bigskip

\centerline{\bf \S 1. Generalities:}

\smallskip

\centerline{\bf $F$--structure vs Poisson structure}

\medskip

{\bf 1.1. Manifolds.} Manifolds considered in this talk
can be $C^{\infty}$, analytic, or formal, eventually with
even and odd coordintes (supermanifolds). The ground field
$K$ of characteristic zero is most often $\bold{C}$ or $\bold{R}$.
Each manifold is endowed with the structure sheaf $\Cal{O}_M$
which is a sheaf of commutative $K$--algebras,
and the tangent sheaf $\Cal{T}_M$ which is a locally
free $\Cal{O}_M$--module of (super)rank equal to the
(super)dimension of $M$. $\Cal{T}_M$ acts on $\Cal{O}_M$
by derivations, and is  a sheaf of Lie (super)algebras
with an intrinsically defined Lie bracket $[\, ,]$. 

\smallskip

There is a classical notion of {\it Poisson structure} on $M$
which endows $\Cal{O_M}$ as well  with a
Lie bracket $\{\, ,\}$ satisfying a certain identity.
Similarly, {\it an $F$--structure} on $M$ endows $\Cal{T}_M$
with an extra operation: (super)commutative and associative
$\Cal{O}_M$--bilinear multiplication $\circ$ with identity $e$. 

\smallskip

In order to describe the structure identities  imposed on these operations,
we recall the notion of the Poisson tensor. Let generally
$A$ be a $K$--linear superspace (or a sheaf of superspaces) endowed with a $K$--bilinear
multiplication and a $K$--bilinear Lie bracket $[\, ,]$. 
Then for any $a,b,c\in A$ put
$$
P_a(b,c):=[a,bc]-[a,b]c-(-1)^{ab}b[a,c] .
\eqno(1.1)
$$
(From here on, $(-1)^{ab}$ and similar notation refers
to the sign occuring in superalgebra when two neighboring elements get
permuted.)

\smallskip

This tensor will be written for $A=(\Cal{O}_M, \cdot, \{\,,\})$
in case of the Poisson structure,  and for $A=(\Cal{T}_M, \circ , [\,,])$
in case of an $F$--structure.

\smallskip

We will now present parallel lists of basic properties of 
Poisson, resp. $F$--manifolds.

\medskip

{\bf 1.2. Poisson (super)manifolds.} (i)${}_P$. {\it Structure identity:}
for all local functions $f,g,h$ on $M$
$$
P_f(g,h)\equiv 0 .
\eqno(1.2)
$$

(ii)${}_P$. Equivalently, each  local
function $f$ on $M$ becomes a local
vector field $X_f$ of the same parity on $M$ via
$X_f(g):=\{f,g\}$

\smallskip

(iii)${}_P$. {\it Maximally nondegenerate case: symplectic structure.}
There exist local canonical coordinates $(q_i,p_i)$ such that
for any $f,g$
$$
\{f,g\}=\sum_{i=1}^n (\partial_{q_i}f\partial_{p_i}g
- \partial_{q_i}g\partial_{p_i}f).
$$
Thus, locally all symplectic manifolds of the same dimension are isomorphic.
Local group of symplectomorphisms is, however,
infinite dimensional.

\medskip

{\bf 1.3. $F$--manifolds.} (i)${}_F$. {\it Structure identity:}
for all local vector fields $X,Y,Z,U$
$$
P_{X\circ Y}(Z,U)=X\circ P_Y(Z,U)+(-1)^{XY}Y\circ P_X(Z,U).
\eqno(1.3)
$$

(ii)${}_F$. Since $(\Cal{T}_M,\circ )$ admits an identity $e$, each  local
vector field on $M$ becomes a local
vector function {\it on the spectral cover $\widetilde{M}$} of $M$ where
$$
\widetilde{M}:=Spec_{\Cal{O}_M}(\Cal{T}_M,\circ )
$$
The last notation is meaningful in each standard category of ringed spaces,
because the tangent sheaf is free as $\Cal{O}_M$--module.
In particular, the structure covering $\widetilde{M}\to M$ is flat.

\smallskip

However, generally $\widetilde{M}$ is not  a manifold
because of nilpotents and/or singularities. Hertling's Theorem 3.6.4
below describes certain important cases when $\widetilde{M}$
is a manifold.

\smallskip

(iii)${}_F$. {\it Maximally nondegenerate case: semisimple $F$--manifolds.}
$\widetilde{M}$  will be a manifold and even an unramified covering
of $M$ in the appropriate ``maximally nondegenerate case'',
namely, when $M$ is pure even, and locally $(\Cal{T}_M,\circ )$ is isomorphic 
to $(\Cal{O}_M^d)$ as algebra, $d=\roman{dim}\,M.$

\smallskip

In this case there exist local canonical coordinates $(u_a)$ (Dubrovin's coordinates)
such that the respective vector fields $\partial_a:=\partial/\partial_a$
are orthogonal idempotents:
$$
\partial_a\circ \partial_a =\delta_{ab}\partial_a .
$$
Thus, locally all semisimple $F$--manifolds of the same dimension are isomorphic.
Local automorphisms of an $F$--semisimple structure
are generated by renumberings and shifts of canonical coordinates:
$$
u_a\mapsto u_{\sigma (a)}+c_a 
$$
so that this structure is more rigid than the symplectic one.

\medskip

{\bf 1.4. Structure embedding of the spectral cover.}  The 
canonical surjective morphism of  sheaves of $\Cal{O}_M$--algebras
$Symm_{\Cal{O}_M}(\Cal{T}_M)\to (\Cal{T}_M,\circ )$ 
induces a closed embedding $\widetilde{M}\to T^*M$.
Its image is Lagrangian (with respect to the canonical
symplectic structure on $T^*M$.) 

\smallskip

As a partial converse, an embedded submanifold $N\subset T^*M$ is
the spectral cover of some generically semisimple $F$--structure
iff $N$ is Lagrangian. 

\medskip

{\bf 1.5. Local decomposition theorem.} For any point
$x$ of a pure even $F$--manifold $M$, the tangent space $T_xM$
is endowed with the structure of a $K$--algebra.
This $K$--algebra can be represented as a direct sum of local
$K$--algebras. The decomposition is unique in the following sense:
the set of pairwise orthogonal idempotent tangent
vectors determining it is well defined.

\smallskip

C. Hertling has shown that this result extends to a
neighborhood of $x$. 

\smallskip

Generally, define the sum of two $F$--manifolds:
$$
(M_1,\circ_1,e_1) \oplus (M_2,\circ_2,e_2) :=
(M_1\times M_2,\circ_1\boxplus\circ_2,e_1\boxplus e_2)
$$
A manifold is called {\it indecomposable} if it cannot be represented
as a sum in a nontrivial way. 

\medskip

{\bf 1.5.1. Theorem.} {\it Every germ $(M,x)$ of a complex
analytic $F$--manifold 
decomposes into a direct sum
of indecomposable germs such that for each summand,
the tangent algebra at $x$ is a local algebra.

\smallskip

This decomposition is unique
in the following sense: the set of pairwise orthogonal idempotent
vector fields determining it is well defined.} 

\medskip

For a proof, see [He], Theorem 2.11. It is based on the following
reinterpretation of (1.3) as an integrability condition:
for any two local vector fields $X,Y$,
$$
\roman{Lie}_{X\circ Y}(\circ )= X\circ \roman{Lie}_Y(\circ )+
(-1)^{XY}Y\circ \roman{Lie}_X(\circ ).
$$
Here $\circ$ is understood as a tensor, to which one can apply
a Lie derivative.

\smallskip

If $(T_xM,\circ )$ is semisimple, this theorem shows that
$M$ is semisimple in a neighborhood of $x$, and
implies the existence and uniqueness (up to reindexing and
shifts) of Dubrovin's canonical coordinates.

\smallskip

For $F$--manifolds with a compatible flat structure,
there exists a considerably more sophisticated operation
of {\it tensor product.} This will be explained in the last section.

\bigskip

\centerline{\bf \S 2. Compatible flat structures}

\smallskip

\centerline{\bf and Euler fields}

\medskip

{\bf 2.1. Flat structures.} Let $M$ be a (super)manifold. 
An (affine) flat structure on $M$ is given by any
of the following data:

\smallskip

(i) A torsionless flat connection $\nabla_0:\,\Cal{T}_M\to
\Omega^1_M\otimes_{\Cal{O}_M}\Cal{T}_M $.

\smallskip

(ii) A local system $\Cal{T}_M^f\subset \Cal{T}_M$ of flat vector fields,
which forms a sheaf of supercommutative Lie algebras of rank $\roman{dim}\,M$
such that $\Cal{T}_M=\Cal{O}_M\otimes_K \Cal{T}_M^f$.

\smallskip

(iii) An atlas whose transition functions are affine linear.

\medskip

{\bf 2.2. Compatibility of a flat structure with a multiplication.}
Assume that $\Cal{T}_M$ is endowed with an $\Cal{O}_M$--bilinear
(super)commutative
and associative multiplication $\circ$,
eventually with unit $e$. At this stage, we do not assume
that (1.3) holds. 

\medskip

{\bf 2.2.1. Definition.} {\it a) A flat structure $\Cal{T}_M^f$
on $M$ is called compatible with $\circ$, 
if in a neighborhood of any point there exists a vector field
$C$ such that for arbitrary local flat vector fields $X,Y$
we have
$$
X\circ Y= [X,[Y,C]].
\eqno(2.1)
$$
$C$ is called a local vector potential for $\circ$.

\smallskip
b)  $\Cal{T}_M^f$ is called compatible with $(\circ ,e)$, if
a) holds and moreover,  $e$ is flat.}
 
\smallskip

{\bf 2.2.2. Proposition.} {\it If $\circ$ admits a compatible flat structure,
then it satisfies the structure identity (1.3) so that
$(M,\circ )$ is an $F$--manifold.}

\smallskip

For a proof, see [Ma2], Prop. 2.4.

\medskip

{\bf 2.3. Pencils of flat connections.} Consider the following
input data:

\smallskip

(i) A flat structure $\nabla_0:\,\Cal{T}_M\to
\Omega^1_M\otimes_{\Cal{O}_M}\Cal{T}_M$ on $M$.

\smallskip

(ii) an (odd) global section 
$A\in \Omega^1_M\otimes_{\Cal{O}_M} End(\Cal{T}_M)$ (Higgs field).

\medskip

Produce from it the following  output data: 

\smallskip

(iii) A pencil of connections $\nabla_{\lambda}=\nabla_{\lambda}^A$ on $\Cal{T}_M$:
$$
\nabla_{\lambda}:= \nabla_0 + \lambda \,A.
$$

\smallskip

(iv) An $\Cal{O}_M$--bilinear composition law 
$\circ =\circ^A$ on $\Cal{T}_M$:
$$
X\circ^A Y:=i_X(A)(Y),\quad i_X(df\otimes G):=Xf\cdot G
$$

\medskip

{\bf 2.3.1. Proposition.} {\it $(M,\circ^A,\nabla_0)$
is an $F$--manifold with compatible flat structure if and
only if $\nabla_{\lambda}^A$ is a pencil
of torsionless flat connections.

\smallskip

In this case, $(M,\circ^A,\nabla_{\lambda}^A)$
is an $F$--manifold with compatible flat structure
for any $\lambda$ as well.}

\medskip

{\bf 2.4. Euler fields.} Consider the following problem:

\smallskip

{\it Given $(M,\circ,\nabla ,e)$ as in 2.2.1 and
the respective pencil of torsionless flat connections
$\nabla_{\lambda}$, extend it to a flat connection 
on $pr_M^*(\Cal{T}_M)$ over $\widehat{M}:=M\times (\lambda$-space).}

\smallskip

The missing data is an extra covariant derivative
in the $\lambda$--direction
$$
\widehat{\nabla}_{\partial_{\lambda}}(Y) = \partial_{\lambda}Y+H(Y)
$$
where $H$
is an even endomorphism of $pr_M^{*}(\Cal{T}_M))$ 
depending on $\lambda$.

\medskip

{\bf 2.4.1. Proposition.} {\it Assume that $\nabla e=0$ and $E:=H(e)$
does not depend on $\lambda$. Then  $\widehat{\nabla}_{\partial_{\lambda}}$
defines a flat extension of $\nabla_{\lambda}$ if and only if $H$
can be given by the formula 
$$
H(X)=X\circ E+ \lambda^{-1} \left(\nabla_XE-X\right),
$$
and $E$ is a vector field on $M$ satisfying the following conditions:
$$
P_E(X,Y)=X\circ Y,\quad [E,\roman{Ker}\,\nabla ] \subset \roman{Ker}\,\nabla . 
$$
}

\smallskip

This means that $E$ is an Euler field of weight one in the following sense:

\medskip

{\bf 2.4.2. Definition.} {\it An Euler field $E$ of weight $d_0$
for $(M,\circ,\nabla ,e)$ is defined by the conditions: $E$ is even, and
$$
P_E(X,Y)=d_0X\circ Y,\quad [E,\roman{Ker}\,\nabla ] \subset \roman{Ker}\,\nabla
$$
}
\smallskip

One easily checks that local Euler fields form a sheaf 
of linear spaces and Lie algebras.
Commutator of two Euler fields is an Euler field of weight zero.
Identity is an Euler field of weight zero.

\smallskip

For proofs and further details, see [Ma2].

\bigskip

\centerline{\bf \S 3. Frobenius manifolds, moduli spaces, mirrors}

\medskip

{\bf 3.1. Definition.} {\it A Frobenius manifold is an
$F$--manifold endowed with a compatible flat structure $\nabla$
and a pseudo--Riemannian metric $g:\,S^2(\Cal{T}_M)\to \Cal{O}_M$
such that

\smallskip

(i) $g$ is flat, and $\nabla$ = the Levi--Civita connection of $g$.

\smallskip

(ii) $g(X\circ Y,Z)= g(X, Y\circ Z).$

\smallskip 

An Euler field $E$ is called compatible with the Frobenius structure if

\smallskip

(iii) $Lie_Eg=Dg$ for a constant $D$.}

\medskip

{\bf 3.2. Associativity equations for the $F$- and Frobenius structures.}
Let us start with an $F$--manifold endowed with a compatible flat structure $\nabla$.  If we choose a
local flat coordinate system $(x^a)$ and write the local vector
potential $C$ as
$C=\sum_c C^c\partial_c$, $X=\partial_a,\,Y=\partial_b$, then (2.1) becomes 
$$
\partial_a\circ\partial_b =\sum_c C_{ab}{}^c\partial_c,\quad
C_{ab}{}^c= \partial_a\partial_bC^c . 
\eqno(3.1)
$$
If moreover
$e$ is flat, we may choose local flat
coordinates so that  $e=\partial_0$, and the
conditions $e\circ \partial_b=\partial_b$ reduce to 
$C_{0b}{}^c=\delta_b^c.$

\smallskip

 If we now choose an arbitrary $C$ and {\it define} a composition
$\circ :\,\Cal{T}_M^f\otimes \Cal{T}_M^f\to \Cal{T}_M$ by the formula
(2.1), it will be automatically supercommutative
in view of the Jacobi identity. Associativity, however, is a quadratic differential constraint on $C$ which was called the {\it oriented associativity
equations} in [LoMa2]: for any $a,b,c,f$
$$
\sum_e C_{ab}{}^eC_{ec}{}^f=(-1)^{a(b+c)} \sum_e C_{bc}{}^eC_{ea}{}^f .
\eqno(3.2)
$$

\smallskip

The choice of a flat invariant metric $(g_{ab})$ as in 3.1 allows us
to reduce $d$ functions $C^a$ to one potential $\Phi$
such that $C^a =\sum_b \partial_b\Phi g^{ba}$. Associativity
equations for $\Phi$ are also called the WDVV equations:
$$
\forall\ a,b,c,d:\quad
\sum_{ef}\Phi_{abe}g^{ef}\Phi_{fcd}=
(-1)^{a(b+c)}\sum_{ef}\Phi_{bce}g^{ef}\Phi_{fad} .
\eqno(3.3)
$$
Here $\Phi_{abc}:=\partial_a\partial_b\partial_c\Phi .$
The flat identity  $e=\partial_0$ equation reduces to
$\Phi_{0ab}=g_{ab}$. The Euler field equations take the form
$$
 E\Phi=(d_0+D)\Phi+(\le 
\,\roman{quadratic\ terms}),\quad Eg=Dg.
$$

\medskip

{\bf 3.3. Semisimple case.} Assume now that $M$ is endowed with a structure
of semisimple Frobenius manifold. Then the multiplication
$\circ$ takes a simple form in canonical local
coordinates $(u^i), \partial_i:=\partial/\partial_{u^i}$, $\partial_i\circ\partial_j=\delta_{ij}.$ Put
$e_j=\partial_j$ so that $e=\sum e_i.$

\smallskip

Instead, the flat metric becomes an unknown functional variable.
It diagonalizes in the canonical coordinates, and
its choice can be reduced to that of {\it a local metric potential}
$\eta = \eta (u)$ such that
$g(\partial_i,\partial_j)=\delta_{ij}\eta_i,$
$\eta_i :=\partial_i\eta .$

\smallskip

Flatness and existence of $\Phi$ reduce to the Darboux--Egoroff equations on $\eta$:
$$
\gamma_{ij}:=\frac{1}{2}\,\frac{\eta_{ij}}{\sqrt{\eta_i\eta_j}},\quad
e_k\gamma_{ij}=\gamma_{ik}\gamma_{kj},\ e\gamma_{ij}=0.
$$

\smallskip

The flat identity equation becomes
$$
e\eta =\roman{const}\, .
$$

\smallskip

Finally, the Euler field equation takes form
$$
E\eta =(D-d_0)\eta +\roman{const}\, .
$$
\medskip

{\bf 3.4. Examples: quantum cohomology.} The simplest examples of 
(formal) Frobenius manifolds coming from algebraic
geometry are furnished by the  quantum cohomology of projective spaces.

\smallskip

Put $H:=H^*(\bold{P}^r,\bold{C})=\oplus_{a=0}^{r}\bold{C}\Delta_a,$
$\Delta_a$:= dual class of $\bold{P}^{r-a}.$

\smallskip

Denote by $(x^a)$  linear coordinates on $H$ dual to $(\Delta_a).$
Take the Poincar\'e pairing for the flat metric:
$g_{ab}=\delta_{a+b,r}$. Obviously, $x^a$ are flat coordinates.

\smallskip

One can construct the following formal solution to the WDVV equations:
$$
\Phi^{\bold{P}^r}(x):= \frac{1}{6}\,(\sum_a x^a\Delta_a)^3+
\sum_{d,n_a\ge 0}N(d;n_2,\dots ,n_r)\,
\frac{(x^2)^{n_2}\dots (x^r)^{n_r}}{n_2!\dots n_r!}\,e^{dx^1}
$$
Here $N(d;n_2,\dots ,n_r)$ is the number of rational curves 
of degree $d$ in $\bold{P}^r$ intersecting $n_a$ hyperplanes
of codimension $a$ in general position. Purely formally, one can prove the following statement:

\medskip

{\bf 3.4.1. Proposition.} {\it (i) There exists a unique solution 
to the WDVV equations of this form
satisfying  $N(1;0,\dots ,0,2)=1.$

\smallskip

(ii) $e:=\partial/\partial x^0$ is the flat identity.

\smallskip

(iii) $E=\sum_a (1-a)x^a\partial_a+(r+1)\partial_1$ is an Euler field
of weight one, $Eg=(2-r)g.$}

\medskip

In much greater generality, the theory of Gromov--Witten invariants endows 
the cohomology space of any smooth projective or compact
symplectic manifold with the canonical  structure
of formal Frobenius manifold. For algebro--geometric aspects of this construction,
see [Ma1].

\smallskip

Recently A. Givental defined and developed quantum $K$--theory:
cf. [Gi] and [Lee]. His construction produces $F$--manifolds 
with a flat invariant metric which are not quite
Frobenius because $e$ is not flat: see [Lee]. Dubrovin's 
almost Frobenius manifolds ([Du2], sec. 3, Definition 9)
are $F$--manifolds with nonflat identity as well.
  
\medskip

{\bf 3.5. Examples: unfolding spaces of isolated singularities.} 
Another important class of examples comes from K. Saito's
theory of unfolding singularities. Again, the simplest
example, that of singularity $z^{n+1}$ at $z=0$, can be described directly and explicitly.

\smallskip

The unfolding space can be identified with the space of
polynomials
$$
 \{p(z)=z^{n+1}+a_1z^{n-1}+\dots +a_n\} = \{(a_i)\in \bold{C}^n\}.
$$
Consider the following space $M$: a point of $M$ is an ordering
of the roots $(\rho_i)$ of $p^{\prime}(z)=0$, $i=1,\dots ,n$,
supplemented by the value of $a_n$. 
(We exclude the case of multiple roots.) The functions
$u^i:=p(\rho_i)$ are local coordinates on $M$. 

\medskip

{\bf 3.5.1. Proposition.} {\it Put $e_i=\partial/\partial u^i$
and define $\circ$ by $e_i\circ e_j=\delta_{ij}e_i.$

\smallskip

Put 
$$
g=\sum_i \frac{(du^i)^2}{p^{\prime\prime}(\rho_i)}\,.
$$
Then $(M, g, \circ )$ is a semisimple Frobenius manifold with metric potential
$$
\eta =-\frac{1}{2(n-1)}\,\sum_i\rho_i^2\,,
$$
flat identity $e =\sum e_i =\frac{\partial}{\partial a_n} ,$
and Euler field 
$$
E = \sum u^ie_i=\dfrac{1}{n+1}
\sum_{i=1}^n (i+1)a_i\dfrac{\partial}{\partial a_i}\, .
$$}

\medskip

K.~Saito's theory endows any unfolding space of an
isolated hypersurface singularity with a canonical structure of
an irreducible generically semisimple $F$--manifold with an Euler field, admitting many structures of a Frobenius manifold which are determined
by a choice of the so called {\it primitive form.}
Existence of primitive forms is a dificult
fact. See [He] for a comprehensive exposition of
K. Saito's theory in the context of $F$--manifolds.

\smallskip

The $F$--structure itself, and the construction of an
Euler field, is furnished by an elementary 
construction which can be easily generalized to include
geometric situations essential for  mirror constructions.
Cf. [Ma1], III.8, where this generalization was called
{\it K. Saito's framework.}

\medskip

{\bf 3.6. Mirror isomorphisms.} The simplest manifestations
of the famous mirror phenomenon, discovered by physicists, 
are highly nontrivial
isomorphisms between certain Frobenius manifolds
representing quantum cohomology ($A$--side) and other
Frobenius structures carried by generalized,
or extended, moduli spaces ($B$--side) and their appropriate
subspaces, for example, cut out by the condition
of integrality of the spectrum of the operator $\roman{ad}\,E$
acting upon flat vector fields.
In this subsection, I will briefly review some of the aspects
of this fascinating story directly involving
$F$--manifolds.

\medskip

{\bf 3.6.1. Extended moduli spaces.} 
Roughly speaking, an unobstructed deformation problem
(like that of a Calabi Yau complex structure $X$) classically is governed
by a certain cohomology space $H^d$ of a fixed weight associated to $X$ which
becomes the tangent space to the moduli space at
the point corresponding to $X$.

\smallskip

Any known construction of what can be considered
as an extended moduli space produces the tangent space
which is the full cohomology algebra $H^*$.

\smallskip  

The existing definitions of
extended moduli spaces of, say, complex structures,
are not really satisfactory in at least  two
respects. First,  generally they furnish only a formal
neighborhood of the classical moduli space. Second, it is
not at all clear what kind of structures are parametrized
by their nonclassical points. Any progress in this problem
necessarily envolves enlarging the geometric universe
so as to include dg--spaces (to treat obstructed cases as well), 
noncommutative spaces etc. For a review, see [Me1], \S 4, and  
references therein.

\smallskip

In [Me1], S.~Merkulov works with $\roman{dg}$--manifolds,
and defines for   them the notion of {\it cohomology $F$--structure.}
He then proves the following two results which show the ubiquity of
$F$--structures:

\medskip

{\bf 3.6.2. Theorem.} {\it (i) The smooth part of the extended
moduli space of deformations of complex or symplectic structure
on a compact manifold is endowed with a canonical
$F$--structure.

\smallskip

(ii) If the $\Cal{G}_{\infty}$ operad acts upon a dg vector space
$(V,d)$, then the formal graded manifold which is the completion
of $H^*(V,d)$ has a canonical structure of the cohomological
$F$--manifold.}

\medskip

{\bf 3.6.3. Unfolding spaces vs extended moduli spaces.} 
Let now $f(x_0,\dots ,x_n)$ be a {\it homogeneous polynomial}
with an isolated singularity at zero. Denote by
$X$ the projective hypersurface $f=0,$ i.e. the base
at infinity of the affine cone $f=0.$ An unfolding of $f$
can be written in the form  $f+\sum_i t_ig_i$ where $g_i$
span $\bold{C}[x]$ modulo the Jacobian ideal.
The unfolding contains a subspace consisting of the homogeneous
polynomials of the same degree as $f$. Geometrically,
it parametrizes unfoldings which correspond to the projective deformations
of $X$ and therefore can be associated with at least 
a part of the classical moduli space of deformations of $X$.
Other points deform $f$ and break the affine cone structure of $f=0$.

\smallskip

Thus, this picture is parallel to the one we would like
to see in the extended deformation theory.
From the general perspective of deformation theory, 
it looks somewhat naive and ad hoc.
On the other hand, it does have some desirable
features: the unfolding space is global (or at least,
not just formal in certain directions), and it is
quite clear what geometric objects correspond to the
nonclassical points.

\smallskip

Furthermore, K. Saito's construction provides the $F$--structure
and many Frobenius structures in this context as well.

\smallskip

For applications to the mirror isomorphisms it is thus
important to have an intrinsic characterization of
$F$--manifolds which can be obtained in this way.
At least locally, it is furnished by the following beautiful theorem
due to C.~Hertling ([He], Theorems 5.3 and 5.6).
All spaces below are germs of complex analytic spaces.

\medskip

{\bf 3.6.4. Theorem.} {\it (i) The spectral cover space $\widetilde{M}$
of the $F$--structure on the germ of the unfolding space
of an isolated hypersurface singularity is smooth.

\smallskip

(ii) Conversely, Let $M$ be an irreducible germ of a generically
semisimple $F$--manifold with the smooth spectral cover $\widetilde{M}.$
Then it is (isomorphic to) the germ of the unfolding space
of an isolated hypersurface singularity. Moreover, any isomorphism
of germs of such unfolding spaces compatible with their
$F$--structure comes from a stable right equivalence of the
germs of the respective singularities.}

\smallskip

Recall that the stable right equivalence is generated by adding sums
of squares of coordinates and making invertible analytic
coordinate changes.

\smallskip

In view of this result, it would be interesting to understand
the following

\medskip

{\bf 3.6.5. Problem.} {\it Which quantum cohomology Frobenius spaces have
smooth spectral covers?}

\smallskip

For a discussion of generically simple quantum cohomology
and Dubrovin's conjectures, see [BeMa].

\bigskip

\centerline{\bf \S 4. Pencils of flat connections on an external bundle}

\smallskip

\centerline{\bf and Dubrovin's twisting}

\medskip

{\bf 4.1. Pencils of flat connections on an external bundle.}
In Proposition 2.3.1, we have characterized $F$--manifolds with
a compatible flat connection in terms of pencils
of flat connections on the tangent sheaf.

\smallskip

Let us now weaken the input data: consider a locally free sheaf $\Cal{F}$ on a manifold
$M$ and a pencil (affine line) $\Cal{P}$ of connections
$\nabla:\,\Cal{F}\to
\Omega^1_M\otimes_{\Cal{O}_M}\Cal{F}$. The difference
of any two connections is  an odd global section
$\nabla_1-\nabla_0=A\in \Omega^1_M\otimes_{\Cal{O}_M} End\,\Cal{F}$.
Any two differences are proportional with a constant coefficient; we will often choose
one arbitrarily.
Any $\nabla$ can be extended to an odd derivation 
(denoted again $\nabla$) of $\Omega^*_M\otimes_{\Cal{O_M}}T(\Cal{F})$
in the standard way, where $T(\Cal{F})$ is the total
tensor algebra of $\Cal{F}$.

\smallskip

Assume that all connections in $\Cal{P}$ are flat, that is
for any $\nabla \in \Cal{P}$, $\nabla A=0$ and $A\wedge A =0$.
Again, in view of the De Rham lemma, $A$ can be everywhere
locally written as $\nabla B$ where $B=B_{\nabla}$ (it generally depends
on $\nabla$) is an even section
of $End\,\Cal{F}$. In this setting, we cannot assume $\Cal{P}$
to be torsionless: this makes no sense for an external bundle.
To produce an $F$--structure, we will have instead to fix an
additional piece of data.

\smallskip

Choose a coordinate neighborhood $U$ in $M$ over which
the linear superspace $F$ of local $\nabla$--flat sections of
$\Cal{F}$ trivializes $\Cal{F}.$ Denote by $\widetilde{F}$
the affine supermanifold associated with $F$. 
Denote by $q:\,\widetilde{\Cal{F}}\to M$ be the fibration
of supermanifolds which is ``the total space'' of
$\Cal{F}$ as a vector bundle: for example, in
algebraic geometry this is the relative affine spectrum
of $\roman{Symm}_{\Cal{O}_M}\,(\Cal{F}^*)$, $\Cal{F}^*$
being the dual sheaf. Hence local sections of $\Cal{F}$
become local sections of $q$. Clearly, $\nabla$ trivializes $q$ over $U$: we have
a well defined isomorphism $q^{-1}(U)=\widetilde{F}\times U$
turning $q$ into projection. 

\smallskip

Let now $u\in F$ be a $\nabla$--flat section of $\Cal{F}$
over $U$, $\nabla B=A$, $B\in \roman{End}\,\Cal{F}.$ 
Then $Bu$ is a section of
$\Cal{F}$; we will identify it with a section of $q$
as above. Projecting to $\widetilde{F}$,
we finally get a morphism $Bu:\,U\to \widetilde{F}$.
Thus $Bu$ denotes several different although closely related objects.
If we choose a basis of flat sections in $F$
and a system of local coordinates $(x^a)$ on $M$,
$B$ becomes a local even matrix function $B(x)$ acting from the left
on columns of local functions. Then $Bu$ becomes
the map $U\to \widetilde{F}:\, x\mapsto B(x)u$.
Since $B$ is defined up to $\roman{Ker}\,\nabla$,
this map is defined up to a constant shift.

\medskip

{\bf 4.1.1. Definition.} {\it A section
$u$ of $\Cal{F}$ is called a primitive section
with respect to  $\nabla \in\Cal{P}$ ,
if it is $\nabla$--flat, and $Bu$ is a local
isomorphism of $U$ with a subdomain of
$\widetilde{F}$.}

\smallskip

Since $F$ is linear, $\widetilde{F}$ has a canonical flat structure
$\nabla^F:\,\Cal{T}_{\widetilde{F}}\to \Omega^1_{\widetilde{F}}
\otimes \Cal{T}_{\widetilde{F}}.$ If $u$ is primitive,
$Bu$ is a local isomorphism allowing us to identify locally
$(Bu)^*(\Cal{T}_{\widetilde{F}}) = \Cal{T}_M$. Moreover, we can
consider the pullback of $\nabla^F$ with respect to $Bu$:
$$
\nabla^*:= (Bu)^*(\nabla^F): \Cal{T}_M\to \Omega^1_M\otimes\Cal{T}_M.
\eqno(4.1)
$$ 
Clearly, the local flat structures
induced by the maps $Bu$ on $M$ do not
depend on local choices of $B$ and glue together to
a flat structure on all of $M$ determined by $\nabla$ and $u$.
It follows also that $\nabla^*$ is flat and torsionless.

\smallskip

There is another important isomorphism produced by this
embedding. Namely, restricting $(Bu)^*$ to $F\subset \Cal{T}_{\widetilde{F}}$
we can identify $(Bu)^*(F)\subset \Cal{T}_M$ with $F\subset \Cal{F}$
and then by $\Cal{O}_M$--linearity construct an
isomorphism $\beta^*:\Cal{F}\to\Cal{T}_M$. The connection
$\nabla^*$ identifies with $\nabla$ under this isomorphism.
The inverse isomorphism $(\beta^*)^{-1}$ can be described
as $X\mapsto i_X(\nabla\,Bu)$.

\smallskip

Denote $\Cal{P}^*$ the pencil of flat connections $\beta^*(\Cal{P})$.

\medskip

{\bf 4.1.2. Example.} Assume that $(M,\circ ,e, \nabla_0)$
is an $F$--manifold endowed with a compatible flat
structure and a $\nabla_0$--flat identity. Put $\Cal{F}:=\Cal{T}_M$
and construct $A$ so that $X\circ Y=i_X(A)(Y).$
Then $e$ is a primitive section which induces exactly the initial
flat structure $\nabla_0^*=\nabla_0$. 

\smallskip

In fact, let $(x^a)$ be a $\nabla_0$--flat local
coordinate system on $M$ such that $e=\partial_0$, and 
$(\partial_0, \dots ,\partial_m)$ the dual basis of flat 
vector fields. Then an easy argument shows that as a vector field,
$Be=\sum_a (x^a+c^a)\,\partial_a$ where $c^a$ are constants.

\smallskip

A converse statement holds as well.

\medskip

{\bf 4.2. Theorem.} {\it Let $(M,\Cal{F},\nabla ,A,u)$
be a pencil of flat connections on an
external bundle endowed with a primitive section.
Then $\Cal{P}^*$ is a pencil of torsionless flat connections,
and $e:=\beta^*(u)$ is an identity for one of
the associated $F$--manifold structures $\circ$.}

\smallskip

For a proof, see [LoMa2], \S 5, and [Ma2], \S 4.

\medskip

{\bf 4.3. The tangent bundle considered as an external bundle.}
Let now $(M,\Cal{P})$ be a manifold with a unital pencil
of torsionless flat connections as in sec. 2.3. Assume moreover
that one (hence each) identity $e_A$ is flat with respect to some
$\nabla_0 \in \Cal{P}.$ Choose as origin  $\nabla_0$
and a coordinate $\lambda$ on $\Cal{P}$ so that
$(M,\Cal{P})$ determines an $F$--manifold structure $\circ$
with $\nabla_0$--flat identity $e$ and multiplication tensor $A$. 

\smallskip

It is natural to study
the family of all $F$--manifold structures that can be obtained from
this one by treating $\Cal{T}_M$ as an external bundle with the pencil $\Cal{P}$,
choosing different $\nabla \in \Cal{P}$ and different
$\nabla$--flat primitive sections, and applying Theorem 4.2.

\smallskip

Let us call an even global vector field $\varepsilon$
{\it a virtual identity}, if it is invertible
with respect to $\circ$ and its inverse $u:=\varepsilon^{-1}$
belongs to $Ker\,\nabla$ for some $\nabla\in\Cal{P}.$   

\medskip

{\bf 4.3.1. Proposition.} {\it a) Inverted virtual identities 
in  $Ker\,\nabla$
are exactly primitive sections of $(\Cal{T}_M,\nabla )$
considered as an external bundle. 

\smallskip

b) The map $(\beta^*)^{-1}$ sends a local vector field $X$ to $X\circ u=X\circ\varepsilon^{-1}$.}

\medskip

{\bf 4.4. Dubrovin's twisting.} Let now $(M, \Cal{P},e,\varepsilon )$
be a flat pencil on $\Cal{T}_M$ such that the respective
$F$--manifold is endowed with a
$\nabla_0$--flat identity $e$, and a $\nabla$--flat
inverse  virtual identity $\varepsilon^{-1}$.

\medskip

{\bf 4.4.1. Theorem.} {\it a) Denote  by $*$ the new multiplication
on $\Cal{T}_M$ induced by $(\beta^*)^{-1}$ from $\circ$: 
$$
X*Y=\varepsilon^{-1}\circ X\circ Y.
\eqno(4.2)
$$
Then $(M,*,\varepsilon )$ is an $F$--manifold with identity
$\varepsilon$ (hence our term ``virtual identity'').

\smallskip

b) With the same notation, assume moreover that $\nabla_0\ne \nabla$
and that $A$ is normalized as $A=\nabla -\nabla_0.$ Then
$e$ is an Euler field of weight one for $(M,*,\varepsilon )$.}

\smallskip

For a proof, see [Ma2].

\smallskip

The relationship
between $(M,\circ , e)$ and $(M,*,\varepsilon )$
is almost symmetric, but not quite. To explain this, we
start with a part that admits a  straightforward check:
$$
X\circ Y=e^{*-1}* X*Y
\eqno(4.3)
$$
which is the same as (4.2) with the roles of $(\circ ,e)$
and $(*,\varepsilon )$ reversed. Here $e^{*-1}$ 
denotes the solution $v$ to 
the equation $v*e=\varepsilon$, that is  
$\varepsilon^{-1}\circ v\circ e=
\varepsilon$, therefore $v=e^{*-1}=\varepsilon^{\circ 2}.$
After this remark one sees that (4.3) follows from (4.2).

\smallskip

Slightly more generally, one easily checks that 
the map inverse to
$X\mapsto \varepsilon^{-1}\circ X$ reads
$Y\mapsto e^{*-1}*Y$ as expected.

\smallskip 

If the  symmetry were perfect, we would now expect $\varepsilon$ 
to be an Euler field of
weight one for $(M,\circ ,e)$. However, one can check that this
is not always true.

\smallskip

The reason of this is that the vector field
$e$ is not a virtual identity for $(M,*,\varepsilon )$:
$e^{*-1}=\varepsilon^2$ is not flat with respect to any
$\nabla^*\in \Cal{P}^*$, so that if we start with $(M,*,\varepsilon)$
and construct $\circ$ via (4.3), we cannot apply
Theorem 4.2 anymore.

\smallskip

Dubrovin has introduced this twisting operation
in the context of Frobenius manifolds and called
it (almost) duality in [Du2].

\bigskip

\centerline{\bf \S 5. Toric compactifications and tensor products}

\smallskip

\centerline{\bf  of $F$--manifolds with compatible flat structure}

\medskip

{\bf 5.1. Setup.} Let $M$ be the formal completion of a linear
(super)space $T$ at 0, that is formal spectrum of $\Cal{R}:=\bold{C}[[x^a]].$
Denote by $\nabla$ the flat connection on $M$ with $Ker\,\nabla=T$.

\smallskip
 
As we have already remarked, the classification of
all formal $F$--manifold structures on $M$
compatible with $\nabla$ is equivalent
to the classification  of the solutions
to the matrix differential equation
$$
\nabla C\wedge \nabla C=0,\quad C\in \roman{End}\,T\widehat{\otimes}_{\bold{C}} m,
\eqno(5.1)
$$
where $m$ is the maximal ideal of $\Cal{R}$. Namely, given $C$, put
$$
X\circ Y:=i_X(\nabla C)(Y)=(XC)Y
$$
In this section we will explain the following result from [LoMa2]:

\medskip

{\bf 5.2. Theorem.} {\it (i) Solutions $C$ to (5.1) bijectively correspond to
the representations in $T$ of a certain algebra $H_{*T}$
constructed from homology spaces of permutohedral
toric compactifications.

\smallskip

(ii) This allows one to define a tensor product operation
on formal $F$--manifolds endowed with a compatible flat structure.}

\medskip

The algebra $H_{*T}$ is essentially motivic. Its possible relation
to the motivic fundamental groups (of Hodge--Tate type)
deserves further study.

\medskip

{\bf 5.3. Permutohedral fans.} Let $B$ be a finite
set. An $N+1$--partition $\tau=\{ \tau_1,\dots ,\tau_{N+1}\}$ of $B$,
by definition, is  {\it a totally
ordered set of $N+1$ pairwise disjoint non--empty subsets of $B$ whose union is
$B$.}

\smallskip

 If $N+1\ge 2,$ $\tau$ determines
a well ordered family of $N$ 2--partitions $\sigma^{(a)}$:  
$$
\sigma^{(a)}_1:=\tau_1\cup\dots\cup\tau_{a},\
\sigma^{(a)}_2:=\tau_{a+1}\cup\dots\cup\tau_{N+1},\ a=1,\dots ,N\, .
$$
A sequence of $N$ 2--partitions $(\sigma^{(i)})$ is called
{\it good} if it can be obtained by such a construction.

\smallskip

Put $N_B := \bold{Z}^B/\bold{Z}$, the subgroup being embedded
diagonally. Similarly, $N_B\otimes{\bold{R}}=\bold{R}^B/\bold{R}$.
Vectors in this space (resp. lattice) are functions
$B\to \bold{R}$ (resp. $B\to \bold{Z}$) considered modulo
constant functions. For a subset $\beta\subset B$,  the function $\chi_{\beta}$ =
1 on $\beta$ and 0 elsewhere, determines such a vector.

\smallskip

The fan $\Phi_B$ in $N_B\otimes{\bold{R}}$, by definition, consists of
the following $l$--dimensional cones $C(\tau )$ labeled by all
$(l+1)$--partitions $\tau$ of $B$.

\smallskip

If $\tau$ is the trivial 1--partition, $C(\tau )=\{0\}$.

\smallskip

If $\sigma$ is a 2--partition, $C(\sigma )$ is generated
by $\chi_{\sigma_1}$, or, equivalently, $-\chi_{\sigma_2}$,
modulo constants.

\smallskip

Generally, let $\tau$ be an $(l+1)$--partition, and 
$\sigma^{(i)},\,i=1,\dots,l$,
the respective good family of 2--partitions. Then
$C(\tau )$ is defined as a cone generated by all $C(\sigma^{(i)})$.

\medskip

{\bf 5.4. Permutohedral toric varieties.} Denote by
$L_B$  the compactification
of the torus $(\bold{G}_m)^B/\bold{G}_m$
associated with the fan $\Phi_B$. 
The permutation group of $B$ acts upon it.

\medskip

{\bf 5.5. The algebra $H_*$.}  It is a graded algebra whose
$n$--th component of $H_*$ is the homology
space $H_{*n}:=H_*(L_n)$ where  $L_n$:=
$L_{\{1,\dots ,n\}}$. The homology space
is spanned by  cycles $\mu (\tau )$
indexed, as well as cones of $\Phi_n$, by partitions $\tau$ of
$\{1,\dots ,n\}$.

\smallskip

The multiplication law is given 
in terms of these generators by the following
prescription. Let $\tau^{(1)}$ (resp. $\tau^{(2)}$)
be a partition of $\{1,\dots ,m\}$ (resp. of $\{1,\dots ,n\}$),
then
$$
\mu (\tau^{(1)})\mu (\tau^{(2)}) = \mu (\tau^{(1)}\cup \tau^{(2)})
$$
where the concatenated partition of 
$\{1,\dots ,m,\, m+1,\dots ,m+n\}$ is defined in an obvious way,
shifting all the components of $\tau^{(2)}$ by $m$. 

\medskip

{\bf 5.6. The algebra $H_{*T}$.} As above, let $T$ be a linear superspace.  
Define  $H_*T$  as the 
algebra of symmetric coinvariants of the 
diagonal tensor product
$$
H_{*}T:= \left(\oplus_{n=1}^{\infty} H_{*n}\otimes T^{\otimes n}
\right)_{\bold{S}} .
$$
The symmetric group $\bold{S}_n$ acts upon the $n$--th graded component.

\smallskip

Given a linear representation $\rho:\,H_*T\to \roman{End}\,T$,
{\it its matrix correlators} are defined by
$$
\tau^{(n)}\langle \Delta_{a_1}\dots\Delta_{a_n}\rangle_{\rho} :=
\rho(\mu (\tau^{(n)})\otimes \Delta_{a_1}\otimes\dots\otimes\Delta_{a_n}) .
$$
Here $\tau^{(n)}$ runs over all partitions of $\{1,\dots ,n\}$
whereas $(a_1,\dots ,a_n)$ runs over all maps $\{1,\dots ,n\}\mapsto I:\,
i\to a_i.$

\smallskip

{\it Top matrix
correlators of $\rho$} correspond
to the identical partitions $\varepsilon^{(n)}$ of $\{1,\dots ,n\}$:
$$
\langle \Delta_{a_1}\dots\Delta_{a_n}\rangle_{\rho} :=
\varepsilon^{(n)}\langle \Delta_{a_1}\dots\Delta_{a_n}\rangle_{\rho}\, .
$$
Given ${\rho}$, construct the series
$$
C_{\rho} =
\sum_{n=1}^{\infty}\sum_{(a_1,\dots ,a_n)}
\frac{x^{a_n}\dots x^{a_1}}{n!}\,
\langle \Delta_{a_1} \dots \Delta_{a_n}\rangle_{\rho} \in \bold{C}[[x]]\otimes \roman{End}\,T.
$$ 
\medskip

{\bf 5.7. Theorem.} {\it (i) We have
$$
\nabla C_{\rho}\wedge \nabla C_{\rho}=0.
$$

\medskip

(ii) Conversely, let $\Delta (a_1,\dots ,a_n)\in\roman{End}\,T$ be a
family of linear operators defined for all $n\ge 1$ and all
maps $\{1,\dots ,n\}\to I:\,i\mapsto a_i$. Assume that the parity
of $\Delta (a_1,\dots ,a_n)$ coincides with the sum of the
parities of $\Delta_{a_i}$ and that 
$\Delta (a_1,\dots ,a_n)$ is (super)symmetric with respect
to permutations of $a_i$'s.
Finally, assume that the formal series 
$$
C =
\sum_{n=1}^{\infty}\sum_{(a_1,\dots ,a_n)}
\frac{x^{a_n}\dots x^{a_1}}{n!}\,
\Delta ({a_1}, \dots ,{a_n}) \in \bold{C}[[x]]\otimes \roman{End}\,T
$$
satisfies $\nabla C\wedge \nabla C=0.$ . Then there exists a well defined
representation ${\rho}:\,H_*T\to \roman{End}\,T$ such that
$\Delta ({a_1}, \dots ,{a_n})$ are the top correlators
$\langle \Delta_{a_1} \dots \Delta_{a_n}\rangle_{\rho}$ of this
representation.}

\medskip

{\bf 5.8. Comultiplication and tensor product.} Since components
of $H_*$ are homology groups of compact manifolds, they admit
a natural comultiplication. This allows one to define ring homomorphisms
$$
\Delta_{T_1,T_2}: H_{*T_1\otimes T_2}\to H_{*T_1}\otimes H_{*T_2}.
$$
Hence, if $\rho_i$ is a representation of $H_{*T_i},\ i=1,2$,
the tensor product  $\rho_1\otimes\rho_2$ induces a representation of
$H_{*T_1\otimes T_2}$.

\smallskip

Using Theorem 5.7, we can translate this operation into
a tensor product of formal $F$--manifolds with a compatible
flat structure.

\bigskip

\centerline{\bf References}

\medskip

[BeMa] A.~Bayer, Yu.~Manin. {\it (Semi)simple exercises in quantum cohomology.}
In: The Fano Conference Proceedings,
ed. by A. Collino, A. Conte, M. Marchisio, Universit\'a
di Torino, 2004, 143--173. Preprint math.AG/0103164

\smallskip

[Du1] B.~Dubrovin. {\it Geometry of 2D topological field theory.}
In: Springer Lecture Notes in Math. 1620 (1996), 120--348.

\smallskip

[Du2] B.~Dubrovin. {\it On almost duality for Frobenius manifolds.}
In: Geometry, Topology, and Mathematical Physics. Ed. by V.~Buchstaber,
I. Krichever. AMS Translations, ser. 2, vol. 212.
Providence, Rhode Island, 2004, 75--132. Preprint math.DG/0307374.

\smallskip

[Gi] A.~Givental. {\it On the WDVV--equation in quantum $K$--theory.}
Fulton's Festschrift, Mich. Math. Journ., 48 (2000), 295--304.

\smallskip

[He] C.~Hertling. {\it Frobenius manifolds and moduli spaces for singularities.}
Cambridge University Press, 2002.

\smallskip

[HeMa] C.~Hertling, Yu.~Manin. {\it Weak Frobenius manifolds.}
Int. Math. Res. Notices, 6 (1999), 277--286. Preprint
math.QA/9810132.

\smallskip

[Lee] Y.~P.~Lee. {\it Quantum $K$--theory I: foundations.}
Preprint math.AG/0105014.

\smallskip

[LoMa1] A.~Losev, Yu.~Manin. {\it New moduli spaces of pointed curves and pencils of flat connections.} Fulton's Festschrift,
Michigan Journ. of Math., 48 (2000), 443--472. Preprint math.AG/0001003

\smallskip

[LoMa2] A.~Losev, Yu.~Manin. {\it Extended modular operad.} 
In: Frobenius Manifolds,
ed. by C. Hertling and M. Marcolli,  Vieweg \& Sohn Verlag,
Wiesbaden, 2004, 181--211. Preprint 
math.AG/0301003.

\smallskip

[Ma1] Yu.~Manin. {\it Frobenius manifolds, quantum cohomology, and moduli spaces.} AMS Colloquium Publ. 47, Providence RI, 1999, 303 pp.

\smallskip

[Ma2] Yu.~Manin. {\it $F$--manifolds with flat structure and Dubrovin's duality.}
Preprint math.DG/0402451

\smallskip

[Me1] S.~Merkulov. {\it Operads, deformation theory and $F$--manifolds.}
In: Frobenius Manifolds,
ed. by C. Hertling and M. Marcolli,  Vieweg \& Sohn Verlag,
Wiesbaden, 2004, 213--251. Preprint math.AG/0210478.

\smallskip

[Me2] S.~Merkulov. {\it PROP profile of Poisson geometry.}
Preprint math.AG/0401034.

\enddocument